\documentclass[a4paper,12pt,oneside]{article} 

\title{Generic Vanishing Fails for Surfaces in Positive Characteristic}
\author{Stefano Filipazzi\thanks{The author was partially supported by DMS-1300750, DMS-1265285 and a grant from the Simons Foundation, Award Number 256202.}}
\date{\vspace{-5ex}} 

\newcommand{\Addresses}{{
  \bigskip
  \footnotesize

  \textsc{Department of Mathematics, University of Utah, 155 South 1400 East,
    Salt Lake City, UT 84112-0090, USA}\par\nopagebreak
  \textit{E-mail address}: \texttt{filipazz@math.utah.edu}

}}

\usepackage{amsthm}
\usepackage{amssymb}
\usepackage{amsmath,amscd}
\usepackage[english,italian]{babel}
\usepackage{latexsym}
\usepackage{graphicx,subfigure}
\usepackage[a4paper,top=30mm,bottom=30mm,left=30mm,right=30mm]{geometry}
\usepackage{tikz}
\usetikzlibrary{matrix,arrows,decorations.pathmorphing}

\usepackage{hyperref}

\usepackage{mathtools}

\usepackage[babel]{csquotes}
\usepackage[maxnames=99,style=alphabetic,backend=bibtex]{biblatex}
\renewbibmacro{in:}{%
  \ifentrytype{article}{}{\printtext{\bibstring{in}\intitlepunct}}} 
  
\renewbibmacro*{issue+date}{%
  \ifboolexpr{not test {\iffieldundef{year}} or not test {\iffieldundef{issue}}}
    {\printtext[parens]{%
       \iffieldundef{issue}
         {\usebibmacro{date}}
         {\printfield{issue}%
          \setunit*{\addspace}%
          \usebibmacro{date}}}}
    {}%
  \newunit}

\theoremstyle{definition}
\newtheorem{teo}{Theorem}[section]
\newtheorem{cor}[teo]{Corollary}
\newtheorem{lemma}[teo]{Lemma}
\newtheorem{prop}[teo]{Proposition}
\newtheorem{oss}[teo]{Observation}

\newtheorem*{teo*}{Theorem}
\newtheorem*{def*}{Definition}
\newtheorem*{oss*}{Observation}

\newtheorem*{ex*}{Example}
\newtheorem*{cor*}{Corollary}

\newtheorem*{ack*}{Acknowledgements}

\newtheorem*{dimo1*}{Proof}

\newenvironment{claim}[1]{\par\noindent\underline{Claim:}\space#1}{}
\newenvironment{claimproof}[1]{\par\noindent\underline{Proof:}\space#1}{\hfill $\blacksquare$}

\newcommand{\F}{\mathcal{F}}
\newcommand{\G}{\mathcal{G}}
\newcommand{\Ox}{\mathcal{O}}

\newcommand{\Pro}{\mathbb{P}}

\newcommand{\du}{\,\check{\vrule height1.3ex width0pt}}

\newcommand{\Spec}{\mathrm{Spec}}

\DeclareMathOperator{\codim} {codim}
\DeclareMathOperator{\car} {char}

\begin{document}

\selectlanguage{english}

\maketitle

\begin{abstract}
We show that there exist smooth surfaces violating Generic Vanishing in any characteristic $p \geq 3$. As a corollary, we recover a result of Hacon and Kov\'{a}cs, producing counterexamples to Generic Vanishing in dimension 3 and higher.
\end{abstract}

\section{Introduction}

Vanishing theorems are some of the most powerful tools in the study of algebraic varieties. In particular, Green and Lazarsfeld showed the following fundamental result.
\begin{teo}[Generic Vanishing, \cite{GL1} and \cite{GL2}] \label{GL}
Let $X$ be a smooth complex projective variety. Then every irreducible component of
\begin{equation}
V^i(\omega_X)=\lbrace P \in \mathrm{Pic}^0(X) \,|\, h^i(X,\omega_X \otimes P) \neq 0 \rbrace
\end{equation}
is a translate of a subtorus of $\mathrm{Pic}^0(X)$ of codimension at least
\begin{equation}
i-(\dim(X)-\dim(a_X(X))),
\end{equation}
where $a_X$ denotes the Albanese morphism. If $\dim(X)=\dim(a_X(X))$, then there are inclusions
\begin{equation}
V^0(\omega_X)\supset V^1(\omega_X) \supset \ldots \supset V^{\dim(X)}(\omega_X)=\lbrace \Ox_X \rbrace.
\end{equation}
\end{teo}
This result has played an important role in the classification of irregular varieties, i.e. varieties carrying a non-trivial morphism to an abelian variety. In particular, Generic Vanishing has been heavily used in the classification of irregular varieties in characteristic 0 (see for example \cite{HP02}, \cite{PP09}, \cite{JLT12}).

A natural question is whether an analog of this result holds in positive characteristic, and whether a similar strategy could be effective for the classification of irregular varieties in this setting. First, we notice that, in the case of maximal Albanese dimension, Generic Vanishing admits a formulation in the language of derived categories, due to Hacon \cite{Hac04}, with refinements by Pareschi and Popa \cite{PP11}.
\begin{teo}[Generic Vanishing, \cite{Hac04} and \cite{PP11}] \label{GVHac}
Let $\F$ be a coherent sheaf on an abelian variety $A$ defined over an algebraically closed field $\mathbb{K}$ of arbitrary characteristic. The following are equivalent:
\begin{itemize}
\item $\codim_{\widehat{A}}  V^i(\F)\geq i$ for all $i \geq 0$, where $V^i(\F)=\lbrace y \in \widehat{A}| h^i(A,\F \otimes \mathcal{L}_y)>0 \rbrace$;
\item for any sufficiently ample line bundle $L$ on $\widehat{A}$,
\begin{equation}
H^i(A,\F \otimes \widehat{L\check{\vrule height1.3ex width0pt}}\,)=0 \quad \forall i > 0;
\end{equation}
\item There is an isomorphism
\begin{equation}
Rp_{\widehat{A},*}(p_A^* D_A(\F)\otimes \mathcal{L}) \simeq R^0p_{\widehat{A},*}(p_A^* D_A(\F)\otimes \mathcal{L}).
\end{equation}
\end{itemize}
Here $\mathcal{L}$ denotes the Poincar\'{e} line bundle on $A \times \widehat{A}$, $p_A$ and $p_{\widehat{A}}$ denote the two projections, and $D_A(-)$ denotes the dualizing functor on $D^b(A)$.
\end{teo}

A coherent sheaf on an abelian variety satisfying one of the equivalent conditions in Theorem \ref{GVHac} is said to be a GV-sheaf.

\begin{oss}
Under the assumption that $X$ is smooth projective over a field of characteristic 0, via Theorem \ref{GVHac} we can recover a generalization of Theorem \ref{GL}. The strategy is to apply Koll\'{a}r Vanishing to the decomposition $Ra_{X,*} \omega_X = \bigoplus R^i a_{X,*} \omega_X [-i]$, in order to deduce that $H^i(A,R^ja_{X,*} \omega_X \otimes \widehat{L\check{\vrule height1.3ex width0pt}}\,)=0$ for any $i>0$ and hence each $R^ja_{X,*} \omega_X$ is a GV-sheaf. Also, it is worthwhile to point out that a statement in the flavor of Theorem \ref{GL} can be recovered in positive characteristic under certain stronger assumptions. For a more detailed discussion about this topic, see \cite{Hac04} and \cite{PP11}.
\end{oss}

From now on, when we refer to Generic Vanishing, we will mean Theorem \ref{GVHac}. For our purposes, we should think of the coherent sheaf $\F$ in the statement as $\lambda_* \omega_X$, where $\lambda: X \rightarrow A$ is a generically finite morphism.

In a recent paper, Hacon and Kov\'{a}cs show that Generic Vanishing does not extend to singular varieties, nor to positive characteristic \cite{HK13}. As for the positive characteristic case, their strategy produces counterexamples in dimension at least 3. This is due to both the geometric construction they consider, and to the fact that they rely on the failure of Grauert-Riemenschneider Vanishing. Since this vanishing holds for smooth surfaces in positive characteristic \cite[Theorem 10.4]{Kol13}, the search for a counterexample in dimension 2 has to involve a different strategy.

In the following, our varieties will be defined over an algebraically closed field $\mathbb{K}$ of 
positive characteristic $p$. In this work, we prove the following.
\begin{teo}[Main Result]\label{main result}
For any prime $p \geq 3$ there exists a smooth surface $S$ and a principally polarized abelian surface $(A,\Theta)$ such that
\begin{itemize}
\item there is a finite map $a:S \rightarrow A$ of degree coprime with $p$;
\item there is an ample and effective divisor $H$ on $S$ such that $H^1(S,\Ox_S(-H))\neq 0$;
\item $a_* \omega_S$ is not a GV-sheaf.
\end{itemize}
\end{teo}

\begin{oss}
The requirement $p \neq 2$ is merely technical, as it will be evident in the course of the proof. Indeed, we expect the result should extend without difficulty to $p=2$.
\end{oss}

Since the map $a$ in Theorem \ref{main result} is finite, we recover the failure of an equivalent of Theorem \ref{GL} as well.
\begin{cor}\label{main cor}
Let $S$ and $A$ be as in Theorem \ref{main result}. Consider the morphism 
\begin{equation}
\lambda  \coloneqq  a \times id : S \times \widehat{A} \rightarrow A \times \widehat{A}.
\end{equation}
Let $\mathcal{L}$ denote the Poincar\'{e} line bundle on $A \times \widehat{A}$, and $p_{\widehat{A}}$ the projection 
\begin{equation}
p_{\widehat{A}}:S \times \widehat{A} \rightarrow \widehat{A}.
\end{equation}
Then $Rp_{\widehat{A},*}(\mathcal{P})$ is not a sheaf, i.e.
\begin{equation}
Rp_{\widehat{A},*}(\mathcal{P}) \not\simeq R^2p_{\widehat{A},*}(\mathcal{P}),
\end{equation}
where $\mathcal{P}=\lambda^* \mathcal{L}$.
\end{cor}

The proof of Theorem \ref{main result} relies on two ingredients. The first one is the intuition that failure of Kodaira Vanishing should imply failure of Generic Vanishing. Therefore, following the constructions due to Raynaud and Mukai \cite{Muk13}, we produce a smooth surface $S$ violating Kodaira Vanishing and having a finite morphism $a: S \rightarrow A$ to a principally polarized abelian variety $(A,\Theta)$.

The second ingredient is the aforementioned categorical formulation of Generic Vanishing. This allows us to look for a contradiction by focusing on the groups $H^1(A,a_* \omega_X \otimes \widehat{L\check{\vrule height1.3ex width0pt}}\,)$. In particular, we identify $A$ with $\widehat{A}$ through $\Theta$, and study what happens for $L=\Ox_A(n\Theta)$ and $n$ large. 

Finally, we recover the result by Hacon and Kov\'{a}cs \cite{HK13}. Indeed, considering products of surfaces violating Generic Vanishing and abelian varieties, we get the following.
\begin{cor}\label{corollary}
For any $n\geq 2$ and prime $p\geq3$, there exist a smooth $n$-fold $X$ and an abelian variety $Y$ of the same dimension defined over a field of characteristic $p$ such that
\begin{itemize}
\item $X$ admits a finite map $a:X \rightarrow Y$ of degree coprime with $p$;
\item $a_* \omega_X$ is not a GV-sheaf;
\item $Rp_{\widehat{Y},*}(\mathcal{P}) \not\simeq R^np_{\widehat{Y},*}(\mathcal{P})$, where the notation is as in Corollary \ref{main cor}.
\end{itemize}
\end{cor}

\begin{ack*}
The author would like to thank his advisor Christopher Hacon for suggesting the problem, for his insightful suggestions and the encouragement. He would also like to thank Karl Schwede for the helpful conversations, and Andrew Bydlon for the many times he listened to his doubts and ideas. Finally, he would like to thank Hanna Astephan for the continuous feedback about his writing.
\end{ack*}

\section{A counterexample to KV of maximal Albanese dimension}

It is known that, for any prime characteristic $p>0$, there is a smooth surface $\tilde{X}$ defined over an algebraically closed field of characteristic $p$ that violates KV (Kodaira Vanishing). In particular, $\tilde{X}$ carries an ample and effective divisor $\tilde{D}$ such that $H^1(\tilde{X},\Ox_{\tilde{X}}(-\tilde{D}))\neq 0$. For the explicit construction we refer to \cite{Muk13}. Now, in order to discuss failure of GV (Generic Vanishing), we are interested in such a surface also having maximal Albanese dimension. We will obtain it by base change.

\begin{prop}
For any prime $p \geq 3$ there exists a smooth surface $S$ of maximal Albanese dimension that violates KV. In particular, there are an abelian surface $A$ and a finite map $a: S \rightarrow A$ of degree coprime with $p$.
\end{prop}
\begin{dimo1*}
As explained in \cite[Construction 2.1]{Muk13}, the surface $\tilde{X}$ comes with a morphism $g: \tilde{X} \rightarrow X$ to a Tango curve. Such a map admits a section $F_{\infty}$. Furthermore, the interesting ample class $\tilde{D}$ is given by 
\begin{equation}
\tilde{D}=(k-1)F_{\infty}+g^*D'.
\end{equation}
Here $D=kD'$ is a particular effective divisor on $X$, and $k \geq 2$ is coprime with $p$.
\smallskip
\begin{claim}
The surface $\tilde{X}$ admits a finite map to $\Pro^2$ of degree coprime with $p$.
\end{claim}
\smallskip
\begin{claimproof}
It is enough to show $\tilde{X}$ has an ample class of degree not divisible by $p$. Given a point $Q$ on $X$ and a natural number $m$, we consider the self-intersection
\begin{equation}
(\tilde{D}+g^*(mQ))^2=\tilde{D}^2+2\tilde{D}\cdot g^*(mQ)=\tilde{D}^2+2m(k-1).
\end{equation}
If $\tilde{D}^2$ is coprime with $p$, we are done choosing $m=0$. Otherwise, assuming $p\geq 3$, we can arrange $k-1$ to be not divisible by $p$ either \cite[Example 1.3]{Muk13}. Thus, for a suitable choice of $m$, we produce an ample divisor with the desired property. Therefore, for this technical reason, from now on we will assume $p\geq 3$.
\end{claimproof}
\smallskip

This shows that $\tilde{X}$ admits a finite map to $\Pro^2$ of degree not divisible by $p$. 
We also choose an abelian variety $A$ having the same property, namely carrying a finite morphism to $\Pro^2$ of degree coprime with $p$. Then, we consider the fiber product
\begin{center}
\begin{tikzpicture}
\matrix(m)[matrix of math nodes,
row sep=2.6em, column sep=2.8em,
text height=1.5ex, text depth=0.25ex]
{\tilde{X} & S\\
\Pro^2 & A\\};
\path[->,font=\scriptsize,>=angle 90]
(m-1-2) edge node[auto] {$h$} (m-1-1)
edge node[auto] {$a$} (m-2-2)
(m-1-1) edge node[auto] {$f$} (m-2-1)
(m-2-2) edge node[auto] {$g$} (m-2-1);
\end{tikzpicture}
\end{center}

First, we want to argue that $S$ is normal and CM (Cohen-Macaulay). By the fiber product construction, we have that $h_* \Ox_S=f^*g_*\Ox_A$. By \cite[Proposition 5.4, Corollary 5.5]{KM}, $g_*\Ox_A$ is locally free. Hence, $h_*\Ox_S$ is locally free as well. Again by \cite[Proposition 5.4, Corollary 5.5]{KM}, we get that $S$ is CM.

Now, by replacing $g$ with $\varphi \circ g$, where $\varphi \in \mathrm{PGL}(2,\mathbb{K})$ is a generic element, we may assume that the branch loci of $f$ and $g$ meet properly. Therefore, oustide of a finite set in $\Pro^2$, at least one map among $f$ and $g$ is \'{e}tale. By stability of \'{e}tale morphisms under base change, we get that $S$ is regular out of the finite set lying above the intersection of the branch loci. Thus, $S$ is $R_1$. Since $S$ is a CM surface, it is $S_2$; therefore, by Serre's criterion on normality, $S$ is normal.

Now, we will have to show $S$ is regular as well. By the above analysis, we have to consider just a finite number of points $\lbrace P_1, \ldots , P_e \rbrace \subset \Pro^2$. Again, since we are free to compose $g$ with an element of $\mathrm{PGL}(2,\mathbb{K})$, we may also assume that the ramification loci of $f$ and $g$ are regular over the $P_i$'s, and that the differential maps of both $f$ and $g$ have rank 1 above them.

Fix $P=P_i$ for some $i$; then, let $Q \in \tilde{X}$ and $R \in A$ be such that $f(Q)=g(R)=P$ and both $f$ and $g$ ramify at $Q$ and $R$, respectively. Call $O$ the point of $S$ lying above $Q$ and $R$. Since all maps are finite, we may consider an affine neighborhood $\Spec (T)$ of $P$. Let $\Spec (U)$ and $\Spec (V)$ the respective preimages in $\tilde{X}$ and $A$. Thus, the homomorphisms $T \rightarrow U$ and $T \rightarrow V$ are finite. Let $\mathfrak{p}$, $\mathfrak{q}$ and $\mathfrak{r}$ the maximal ideals defining $P$, $Q$ and $R$, respectively.

Given the assumptions just made, we then write the maximal ideals of the above local rings as follows. We have $\mathfrak{p}=(x,y)$, where $\lbrace x=0\rbrace$ and $\lbrace y=0 \rbrace$ are local equations of the branch loci of $f$ and $g$, respectively. Now, we consider the images of these local parameters under the maps at the level of local rings 

\begin{center}
\begin{tikzpicture}
\matrix(m)[matrix of math nodes,
row sep=2.6em, column sep=2.8em,
text height=1.5ex, text depth=0.25ex]
{\Ox_Q & \\
\Ox_P & \Ox_R\\};
\path[->,font=\scriptsize,>=angle 90]
(m-2-1) edge node[auto] {$\psi$} (m-1-1)
(m-2-1) edge node[auto] {$\varphi$} (m-2-2);
\end{tikzpicture}
\end{center}

By construction, we have that around $Q$ and $R$ the ramification divisors map isomorphically to the respective branch divisors through $P$. Therefore, we have
\begin{itemize}
\item $\varphi (y) = cv^k$, where $c$ is a unit in $\Ox_R$, $\varphi(x)=u$, and $\mathfrak{r}=(u,v)$;
\item $\psi (x) = d \alpha^l$, where $d$ is a unit in $\Ox_Q$, $\psi(y)=\beta$, and $\mathfrak{q}=(\alpha,\beta)$.
\end{itemize}

Our goal now is to show that $\Ox_O$ is regular. Since a local ring is regular if and only if its completion is \cite[\href{http://stacks.math.columbia.edu/tag/07NU}{Tag 07NU}]{stacks-project}, we can study the situation at the level of completed rings. Furthermore, Cohen's structure theorem guarantees that, after completion, the above diagram is

\begin{center}
\begin{tikzpicture}
\matrix(m)[matrix of math nodes,
row sep=2.6em, column sep=2.8em,
text height=1.5ex, text depth=0.25ex]
{\mathbb{K}[[\alpha,\beta]] & \\
\mathbb{K}[[x,y]] & \mathbb{K}[[u,v]]\\};
\path[->,font=\scriptsize,>=angle 90]
(m-2-1) edge node[auto] {$\psi$} (m-1-1)
(m-2-1) edge node[auto] {$\varphi$} (m-2-2);
\end{tikzpicture}
\end{center}

Since $U$, $V$ and $U\otimes_T V$ are all finite $T$-modules, we can apply \cite[\href{http://stacks.math.columbia.edu/tag/00MA}{Tag 00MA}]{stacks-project} and \cite[\href{http://stacks.math.columbia.edu/tag/07N9}{Tag 07N9}]{stacks-project}. This tells us that the completion of $\Ox_O$ coincides with $\mathbb{K}[[\alpha,\beta]]\otimes_{\mathbb{K}[[x,y]]}\mathbb{K}[[u,v]]$. On the one hand, we know it is a local ring of dimension 2 by geometric reasons; on the other hand, we have the relations $u=d\alpha^l$ and $\beta=cv^k$. Therefore, one sees that $\mathbb{K}[[\alpha,\beta]]\otimes_{\mathbb{K}[[x,y]]}\mathbb{K}[[u,v]]/(\alpha\otimes 1,1 \otimes v)\cong \mathbb{K}\otimes_{\mathbb{K}}\mathbb{K}$, and therefore the ring is regular. This shows that $S$ is a smooth surface.

The above argument does not show that $S$ is irreducible. Since $S \rightarrow \Pro^2$ has degree coprime with $p$, we can find an irreducible component of $S$ dominating $\Pro^2$ with degree coprime with $p$. In the following, we will replace $S$ by such an irreducible component; in particular, this dominates both $\tilde{X}$ and $A$ with degrees not divisible by $p$. For convenience, we will denote the chosen component by $S$ as well.

Now, we are left with showing that $S$ violates KV too. In order to do so, we introduce a technical result that will be used many times henceforth.
\begin{lemma}\label{lemmino}
Let $f:X \rightarrow Y$ be a finite and surjective morphism of $n$-dimensional projective varieties defined over an algebraically closed field $\mathbb{F}$. Assume that $Y$ is normal and that $\car\mathbb{F}$ does not divide $\deg f$. Then, given a line bundle $\G$ on $Y$, the cohomology groups $H^i(Y,\G)$ are direct summands of the cohomology groups $H^i(X,f^*\G)$. In particular, if $H^i(Y,\G) \neq 0$, then $H^i(X,f^*\G) \neq 0$.
\end{lemma}
\begin{dimo1*}
Since $\car \mathbb{F}$ does not divide $\deg f$, then $\Ox_Y$ is a direct summand of $f_* \Ox_X$ \cite[Proposition 5.7]{KM}. This, together with the projection formula, tells us that $\G$ is a direct summand of $f_*f^*\G=f_*\Ox_X \otimes \G$. In particular, we have that $H^i(Y,\G)$ is a direct summand of $H^i(Y,f_*f^*\G)$. Now, since $f$ is finite, we have $H^i(X,f^*\G)=H^i(Y,f_*f^*\G)$, and the claim follows. \hfill $\square$
\end{dimo1*}

By construction, we can apply Lemma \ref{lemmino} to the morphism $h:S \rightarrow \tilde{X}$. Thus, we get that $H^1(S,\Ox_S(-H))\neq 0$, where we write $H=h^*\tilde{D}$. Since $h$ is finite, $H$ is still ample. Therefore, $S$ violates KV. \hfill $\square$
\end{dimo1*}
\section{A counterexample to GV}

Following \cite{Hac04}, we recall an important construction in the setting of abelian varieties; for a complete discussion, see \cite{Muk81}. Given an abelian variety $B$, we denote by $\widehat{B}$ its dual abelian variety, and by $\mathcal{L}$ the normalized Poincar\'{e} line bundle on $B \times \widehat{B}$. A nondegenerate line bundle $M$ (i.e. $\chi(M)\neq 0$) on $\widehat{B}$ induces an isogeny \cite[p. 59, p.131]{Mumford}
\begin{equation}
\begin{split}
\phi_M : \widehat{B} &\rightarrow B\\
x &\to \tau_x^*M\du \otimes M,
\end{split}
\end{equation}
where $\tau_x$ denotes translation by $x$. Then, we set
\begin{equation}
\widehat{M}=R^{i(M)}p_{B,*}(p_{\widehat{B}}^*M \otimes \mathcal{L}),
\end{equation}
where $i(M)$ is the WIT index of $M$\footnote{In the following, we will consider $M=L\du$ with $L$ ample. In this case, we will have $i(M)=\dim B$, and $|\chi(M)|=h^0(L)$.}, and we have \cite[Proposition 3.11]{Muk81}
\begin{equation}
\phi_M^*(\widehat{M})\cong \bigoplus_{|\chi(M)|}M\du.
\end{equation}

\subsection{A first reduction}

So far, for any prime characteristic $p \geq 3$, we have constructed a smooth surface $S$ of maximal Albanese dimension that violates KV. In particular, there is a finite morphism $a:S \rightarrow A$ to an abelian surface whose degree is coprime with $p$. From now on, we will assume that $A$ is principally polarized, and that $\Theta$ is a symmetric principal polarization. In the following, we will also use $\Theta$ to identify $A$ with $\widehat{A}$; namely, $A$ and $\widehat{A}$ will be identified under the isomorphism $\phi_\Theta: A \rightarrow \widehat{A}$. The image of $\Theta$, which we will denote by $\widehat{\Theta}=\phi_\Theta (\Theta)$, is still a symmetric principal polarization.

By Theorem \ref{GVHac}, we have that $a_* \omega_S$ is a GV-sheaf if and only if
\begin{equation}
H^i(S,\omega_S \otimes a^* \widehat{L\du}\, )=0 \quad \forall i > 0
\end{equation}
for any sufficiently ample line bundle $L$ on $\widehat{A}$. Looking for a counterexample, we may assume that $L=\Ox_{\widehat{A}}(n\widehat{\Theta})$. These line bundles are particularly nice, since, after the identification between $A$ and $\widehat{A}$, $\phi_{n\widehat{\Theta}}$ coincides with $\mathbf{n}$, the multiplication by $n$ \cite[p. 60]{Mumford}. This, together with the symmetry of $\Theta$, ensures that $\phi_{n\widehat{\Theta}}^*\Theta=n^2\widehat{\Theta}$ \cite[p. 59]{Mumford}. For the reader's convenience, we will keep the notation for $A$ and its dual $\widehat{A}$ separate. Also, this choice makes so that, from now on, we can denote $\phi_{n\widehat{\Theta}}$ by $\mathbf{n}$ without any ambiguity.

Intuition suggests that failure of KV, i.e. $H^1(S,\omega_S \otimes \Ox_S (H)) \neq 0$, and GV, i.e. $H^i(S,\omega_S \otimes a^*\widehat{L\du}\,)=0$ for all sufficiently ample $L$, should not be compatible. We consider $L=\Ox_{\widehat{A}}(n\widehat{\Theta})$, so that $\mathbf{-n}^*(\widehat{L\du}\,)=\bigoplus_{h^0(n\widehat{\Theta})}\Ox_{\widehat{A}}(n\widehat{\Theta})$ \cite[Theorem 3.11]{Muk81}, and we will compare $H$ and $\Theta$. The following construction goes in this direction.
\smallskip
\begin{claim}
We may assume that $H-a^*\Theta$ is ample, and that there is $F \in |H-a^* \Theta|$ smooth.
\end{claim}
\smallskip
\begin{claimproof}
For $k>>0$, the divisor $H-a^*\left( \frac{1}{k} \Theta \right)$ is ample. Without loss of generality, we may assume $k$ is not divisible by $p$. Then, we have the Cartesian diagram

\begin{center}
\begin{tikzpicture}
\matrix(m)[matrix of math nodes,
row sep=2.6em, column sep=2.8em,
text height=1.5ex, text depth=0.25ex]
{S & S_k & \widehat{S}\\
A & \widehat{A} &\\};
\path[->,font=\scriptsize,>=angle 90]
(m-1-2) edge node[auto] {$\varphi$} (m-1-1)
edge node[auto] {$a_k$} (m-2-2)
(m-1-1) edge node[auto] {$a$} (m-2-1)
(m-1-3) edge node[auto] {$\psi$} (m-1-2)
(m-2-2) edge node[auto] {$\mathbf{k}$} (m-2-1);
\end{tikzpicture}
\end{center}
where $S_k=S\times_A \widehat{A}$ is a smooth surface. By Lemma \ref{lemmino}, the pair $(S_k,\varphi^*H)$ violates KV as well. Furthermore, we have that $\varphi^*\left(H-a^*\left(\frac{1}{k}\Theta\right)\right)=\varphi^*H-a_k^*(k\widehat{\Theta})$ is ample. Therefore, $\varphi ^* H - a_k^* \widehat{\Theta}$ is ample as well. For $c>>0$ and not divisible by $p$, $c(\varphi^*H-a_k^*\widehat{\Theta})$ is very ample. Let $E \in |c(\varphi^*H-a_k^*\widehat{\Theta})|$ be a smooth curve, and denote by $\psi:\widehat{S}\rightarrow S_k$ the $c$-sheeted cover branched over $E$ induced by $\varphi^*H-a_k^*\widehat{\Theta}$ \cite[Chapter 3]{EV}. By construction, $\widehat{S}$ is smooth, and $a_k \circ \psi : \widehat{S} \rightarrow \widehat{A}$ is finite of order not divisible by $p$ \cite[Lemma 3.15]{EV}. Furthermore, $(\varphi \circ \psi)^*H-(a_k \circ \psi)^*\widehat{\Theta}$ is integral, ample and admits a smooth element in its linear series (i.e. the divisor given by the preimage of $E$). Lastly, again by Lemma \ref{lemmino}, $(\varphi \circ \psi)^*H$ violates KV. Therefore, up to replacing $a:S \rightarrow A$ by $a_k \circ \psi : \widehat{S} \rightarrow \widehat{A}$, and switching the role between $A$ and $\widehat{A}$, we may assume that $H-a^*\Theta$ is ample, and admits a smooth element in its linear series.
\end{claimproof}

\subsection{Failure of KV implies failure of GV}

As showed above, we have a polarized surface $(S,H)$ and a principally polarized abelian surface $(A,\Theta)$ satisfying the following:
\begin{itemize}
\item there is a finite surjective projective morphism $a:S \rightarrow A$ whose degree is not divisible by $p$;
\item $H^1(S,\Ox_S(-H))\neq 0$;
\item $H-a^*\Theta$ is ample, and there is $F \in |H-a^*\Theta|$ smooth.
\end{itemize}
Now, we are ready to complete the proof of Theorem \ref{main result}.

\begin{dimo1*}[Theorem \ref{main result}]
For $n>>0$ and not divisible by $p$, we consider the following Cartesian diagram

\begin{center}
\begin{tikzpicture}
\matrix(m)[matrix of math nodes,
row sep=2.6em, column sep=2.8em,
text height=1.5ex, text depth=0.25ex]
{S & S_n\\
A & \widehat{A}\\};
\path[->,font=\scriptsize,>=angle 90]
(m-1-2) edge node[auto] {$\varphi$} (m-1-1)
edge node[auto] {$a_n$} (m-2-2)
(m-1-1) edge node[auto] {$a$} (m-2-1)
(m-2-2) edge node[auto] {$\mathbf{n}$} (m-2-1);
\end{tikzpicture}
\end{center}

Let $\widehat{L\du}$ denote the vector bundle induced by $L\du$, where $L=\Ox_{\widehat{A}}(n\widehat{\Theta})$. By the identification given by $\phi_{\Theta}$, we have $\phi_{L\du}=\mathbf{-n}$. Then, by \cite[Proposition 3.11]{Muk81}, $(\mathbf{-n})^*(\widehat{L\du}\,)=\bigoplus_{h^0(n\widehat{\Theta})}\Ox_{\widehat{A}}(n\widehat{\Theta})$. Given that the polarization $\widehat{\Theta}$ is symmetric, we get $\mathbf{n}^*(\widehat{L\du}\,)=\bigoplus_{h^0(n\widehat{\Theta})}\Ox_{\widehat{A}}(n\widehat{\Theta})$.

Recall that $H-a^*\Theta$ is ample with a smooth element $F$ in its linear series. Since $\varphi$ is \'{e}tale, the same holds for $\varphi^*(H-a^*\Theta)=\varphi^*H - n a_n^*(n\widehat{\Theta})$; call $C=\varphi^{-1}(F)$ the smooth curve in the linear series $|\varphi^*H - n a_n^*(n\widehat{\Theta})|$. Also, since $n$ is coprime with $p$, we have that $\varphi^* H$ violates KV by Lemma \ref{lemmino}.

We now consider the following exact sequence
\begin{equation}
0 \rightarrow \omega_{S_n} \otimes \Ox_{S_n}(n a_n^*(n\widehat{\Theta})) \rightarrow \omega_{S_n} \otimes \Ox_{S_n}(\varphi^*H) 
\rightarrow \omega_C \otimes \Ox_{S_n}(n a_n^*(n\widehat{\Theta})) \rightarrow 0.
\end{equation}
By Serre duality $H^1(C,\omega_C \otimes \Ox_{S_n}(n a_n^*(n\widehat{\Theta})))=0$. Therefore, we get a surjection
\begin{equation}
H^1(S_n,\omega_{S_n} \otimes \Ox_{S_n}(n a_n^*(n\widehat{\Theta}))) \rightarrow H^1(S_n,\omega_{S_n} \otimes \Ox_{S_n}(\varphi^*H)) 
\rightarrow 0.
\end{equation}

Now, assume $a_* \omega_S$ is a GV-sheaf. Our goal is to derive a contradiction. We will achieve it by showing that the trivial group surjects onto $H^1(S_n,\omega_{S_n}\otimes \Ox_{S_n}(\varphi^*H))$, which is non-zero by construction. Going in this direction, we consider the group $H^1(S_n,\omega_{S_n} \otimes \Ox_{S_n}(a_n^*(n\widehat{\Theta})))$ and the auxiliary short exact sequence of sheaves that follows.

The divisors $n a_n^*(n\widehat{\Theta})$ and $a_n^*(n\widehat{\Theta})$ differ by $a_n^*(n(n-1)\widehat{\Theta})$. We may assume that $n$ is large enough (e.g. $n \geq 3$ \cite[p. 163]{Mumford}), so that $n(n-1)\widehat{\Theta}$ is very ample on $\widehat{A}$. Let $F \in |n(n-1)\widehat{\Theta}|$ be a generic smooth element, and $E=a_n^*F \in |a_n^*(n(n-1)\widehat{\Theta})|$. Multiplying by the equation of $E$, we get a short exact sequence
\begin{equation}
0 \rightarrow \omega_{S_n} \otimes \Ox_{S_n}(a_n^*(n\widehat{\Theta})) \rightarrow \omega_{S_n} \otimes \Ox_{S_n}(n a_n^*(n\widehat{\Theta})) 
\rightarrow \omega_{S_n} \otimes \Ox_{E}(n a_n^*(n\widehat{\Theta})) \rightarrow 0.
\end{equation}
To get a surjection
\begin{equation}
H^1(S_n,\omega_{S_n} \otimes \Ox_{S_n}(a_n^*(n\widehat{\Theta}))) \rightarrow H^1(S_n,\omega_{S_n} \otimes \Ox_{S_n}(\varphi^*H)) 
\rightarrow 0,
\end{equation}
it is then enough to show $H^1(E,\omega_{S_n} \otimes \Ox_{E}(n a_n^*(n\widehat{\Theta})))=0$.

\smallskip
\begin{claim}
We may assume $E$ is smooth.
\end{claim}
\smallskip
\begin{claimproof}
By Zariski-Nagata purity \cite[\href{http://stacks.math.columbia.edu/tag/0BMB}{Tag 0BMB}]{stacks-project}, the branch locus $B_{a_n}$ of $a_n$ is purely divisorial. Since $F$ is a general element of a very ample linear series, we may assume that it intersects the branch locus properly \cite[Remark 8.18.1]{Har77}. Furthermore, we may assume that that the ramification locus is smooth above the intersection, and that the rank of the differential of $a_n$ is 1 there. Now, we want to show that under these assumptions $E$ is smooth. Since the map is \'{e}tale away from $B_{a_n}$, we just have to focus on where $F$ and $B_{a_n}$ meet. Let $P \in \widehat{A}$ be such a point, and let $Q \in S_n$ be a point in its preimage.

Let $\mathfrak{p}$ be the maximal ideal of $\Ox_P$. By the above assumptions, we have that $\mathfrak{p}=(x,y)$, where $\lbrace x=0 \rbrace$ is a local equation for $B_{a_n}$, and $\lbrace y=0 \rbrace$ is a local equation for $F$. Now, we can argue as we did before. At the level of local rings, we have
\begin{center}
\begin{tikzpicture}
\matrix(m)[matrix of math nodes,
row sep=2.6em, column sep=2.8em,
text height=1.5ex, text depth=0.25ex]
{\Ox_Q & \\
\Ox_P & \Ox_P / (y)\\};
\path[->,font=\scriptsize,>=angle 90]
(m-2-1) edge node[auto] {$\varphi$} (m-1-1)
(m-2-1) edge node[auto] {$\pi$} (m-2-2);
\end{tikzpicture}
\end{center}
In particular, we know that $\varphi(x)=cv^k$, where $c$ is a unit in $\Ox_Q$, $\varphi(y)=u$, and $\mathfrak{q}=(u,v)$. This in particular shows that a local equation for $E$ at $Q$ is $\lbrace u=0 \rbrace$; since $u$ is a local parameter and $S_n$ is smooth, we have that $E$ is regular at $Q$.
\end{claimproof}
\smallskip

Now, since $E$ is smooth, we are free to use the above adjunction argument again and write
\begin{equation}
\omega_{S_n} \otimes \Ox_{E}(n a_n^*(n\widehat{\Theta}))=\omega_E \otimes \Ox_{S_n}(a_n^*(n\widehat{\Theta})).
\end{equation}
In particular, by Serre duality we get $H^1(E,\omega_{S_n} \otimes \Ox_{E}(n a_n^*(n\widehat{\Theta})))=0$.



Therefore, to derive a contradiction, it is enough to show
\begin{equation}
H^1(S_n,\omega_{S_n} \otimes \Ox_{S_n}(a_n^*(n\widehat{\Theta})))=0.
\end{equation} 

Now, we recall that $L= \Ox_{\widehat{A}}(n\widehat{\Theta})$. Then, by \cite[Proposition 3.11]{Muk81}, we have
\begin{equation}
\begin{split}
\bigoplus_{h^0(L)} \widehat{L\du} &= \phi_{L\du,*}L \\ &= (\mathbf{-n})_*L \\
&= (\mathbf{-n})_*\Ox_{\widehat{A}}(n\widehat{\Theta})\\
&= \mathbf{n}_* \Ox_{\widehat{A}}(n\widehat{\Theta}),
\end{split}
\end{equation}
where the last identity follows from the symmetry of $\widehat{\Theta}$.

Then, if GV holds for $a:S \rightarrow A$, we have
\begin{eqnarray}
H^1(S_n,\omega_{S_n}\otimes a_n^*\Ox_{\widehat{A}}(n\widehat{\Theta})) = & \textrm{$\varphi$ \'{e}tale, $\omega_{S_n}=\varphi^* \omega_S$} \nonumber\\ 
H^1(S_n,\varphi^*\omega_{S}\otimes a_n^*\Ox_{\widehat{A}}(n\widehat{\Theta})) = & \textrm{$a_n$ finite, projection formula} \nonumber \\
H^1(\widehat{A},a_{n,*}\varphi^*\omega_{S}\otimes \Ox_{\widehat{A}}(n\widehat{\Theta})) = & \textrm{flat base change} \nonumber \\
H^1(\widehat{A},\mathbf{n}^*a_*\omega_{S}\otimes \Ox_{\widehat{A}}(n\widehat{\Theta})) = & \textrm{$\mathbf{n}$ finite, projection formula} \\
H^1(A,a_* \omega_S \otimes \mathbf{n}_*\Ox_{\widehat{A}}(n\widehat{\Theta})) = & \nonumber \\
\bigoplus_{h^0(L)} H^1(A,a_*\omega_S \otimes \widehat{L\du}\,)= & 0. \phantom{aaaan finite, projection formula}\nonumber 
\end{eqnarray}


This gives the claimed contradiction. Therefore the morphism $a:S \rightarrow A$ above constructed provides a counterexample to GV in positive characteristic. This concludes the proof of Theorem \ref{main result}. \hfill $\square$
\end{dimo1*}



\subsection{Extension to higher dimensions}

The conterexample to GV for surfaces can be manipulated in order to provide counterexamples in dimension 3 and higher. This leads to the proof of the first two parts of Corollary \ref{corollary}.

\begin{dimo1*}[Corollary \ref{corollary}, part 1 and 2]
In the following, fix a smooth surface $S$ as in Theorem \ref{main result}. In particular, there is a finite morphism $a : S \rightarrow A$ to a principally polarized abelian variety $(A,\Theta)$ such that $H^1(A,a_* \omega_S \otimes \widehat{L\check{\vrule height1.3ex width0pt}}\,)\neq 0$. Here we have $L=\Ox_{\widehat{A}}(n\widehat{\Theta})$, where $n>>0$ is not divisible by $p$.

Now, consider an abelian variety $B$ of dimension $m$. We consider the finite morphism
\begin{equation}
f \coloneqq  a \times id_B : S \times B \rightarrow A \times B.
\end{equation}
First, we notice that $\omega_{S \times B}= \omega_S \boxtimes \Ox_B$. Therefore
\begin{equation}
f_* \omega_{S \times B}=a_* \omega_S \boxtimes \Ox_B.
\end{equation}
Now, consider $M=L\boxtimes N$, where $L$ is as above and $N$ is sufficiently ample on $B$. We can assume that $L$ is ample enough, so that Theorem \ref{GVHac} applies to $M$.

Now, we know that the Poincar\'{e} line bundle on $A\times B \times \widehat{A} \times \widehat{B}$ is the box product of the Poincar\'{e} line bundles on $A \times \widehat{A}$ and $B\times \widehat{B}$. Also, we know that Fourier-Mukai transforms behave well with respect to the box product \cite[Exercise 5.13]{Huy}. Therefore, we have
\begin{equation}
\widehat{L\check{\vrule height1.3ex width0pt} \boxtimes M\du}\cong \widehat{L\check{\vrule height1.3ex width0pt}} \boxtimes \widehat{M\du}.
\end{equation}

Hence, by the K\"{u}nneth formula \cite[Proposition 9.2.4]{Kem}, we get
\begin{equation}
\begin{split}
H^1(A \times B, f_* (\omega_{S \times B}) \otimes \widehat{L\check{\vrule height1.3ex width0pt} \boxtimes M\du}\,) &= H^1(A \times B, (a_*\omega_S \otimes \widehat{L\check{\vrule height1.3ex width0pt}}\,) \boxtimes ( \widehat{M\du}\, ))\\
&= \bigoplus_{p+q=1}H^p(A,a_*\omega_S \otimes \widehat{L\check{\vrule height1.3ex width0pt}}\,) \otimes_k H^q(B,\widehat{M\du}\, ).
\end{split}
\end{equation}

Since $H^0(B,\widehat{M\du}\, )$ does not vanish by construction \cite[Proposition 3.11]{Muk81}, and we have $H^1(A,a_*\omega_S \otimes \widehat{L\check{\vrule height1.3ex width0pt}}\,) \neq 0$ by choice of $S$, we conclude that GV fails for $S \times B$. In particular, we have exhibited a smooth $(m+2)$-fold for which GV fails. This proves the first two parts of Corollary \ref{corollary}. \hfill $\square$
\end{dimo1*}

\subsection[Proof of Corollary 1.7]{$Rp_{\protect\widehat{Y},*}(\protect\mathcal{P})$ is not a sheaf}

Finally, we prove Corollary \ref{main cor} and the last part of Corollary \ref{corollary}. Since the former is a particular case of the latter, we will work in the setting of Corollary \ref{corollary}.

\begin{dimo1*}[Corollary \ref{corollary}, part 3]
By the second part of Corollary \ref{corollary}, we know that
\begin{equation}
Rq_{\widehat{Y},*}(q_Y^* D_Y(a_* \omega_X)\otimes \mathcal{L}) \not \simeq R^0q_{\widehat{Y},*}(q_Y^* D_Y(a_* \omega_X)\otimes \mathcal{L}),
\end{equation}
where $\mathcal{L}$ denotes the Poincar\'{e} line bundle on $Y \times \widehat{Y}$, $q_Y$ and $q_{\widehat{Y}}$ denote the two projections, and $D_Y(-)$ denotes the dualizing functor on $D^b(Y)$. Now, consider the morphism
\begin{equation}
\lambda  \coloneqq  a \times id : X \times \widehat{Y} \rightarrow Y \times \widehat{Y}.
\end{equation}
Let $p_{\widehat{Y}}$ be the projection
\begin{equation}
p_{\widehat{Y}}:X \times \widehat{Y} \rightarrow \widehat{Y},
\end{equation}
and set $\mathcal{P}=\lambda^*\mathcal{L}$.

To conclude the proof is then enough to show
\begin{equation}
Rq_{\widehat{Y},*}(q_Y^* D_Y(a_* \omega_X)\otimes \mathcal{L}) \simeq Rp_{\widehat{Y},*}(\mathcal{P})[n].
\end{equation}

Now, we consider the following natural isomorphisms
\begin{eqnarray}
Rq_{\widehat{Y},*}(q_Y^* D_Y(a_* \omega_X)\otimes \mathcal{L}) \simeq & \textrm{$a$ finite} \nonumber\\ 
Rq_{\widehat{Y},*}(q_Y^* D_Y(Ra_* \omega_X)\otimes \mathcal{L}) \simeq & \textrm{Grothendieck-Verdier duality} \nonumber \\
Rq_{\widehat{Y},*}(q_Y^* Ra_*(D_X \omega_X)\otimes \mathcal{L}) \simeq & \textrm{dualizing functor} \nonumber \\
Rq_{\widehat{Y},*}(q_Y^* Ra_*\Ox_X \otimes \mathcal{L})[n] \simeq & \textrm{flat base change} \nonumber \\
Rq_{\widehat{Y},*}(R\lambda_* p_X^*\Ox_X \otimes \mathcal{L})[n] \simeq & \textrm{projection formula, $\lambda$ finite}\\
Rq_{\widehat{Y},*}(R\lambda_*( p_X^*\Ox_X \otimes \lambda^* \mathcal{L}))[n] \simeq & \nonumber \\
Rq_{\widehat{Y},*}(R\lambda_*\mathcal{P})[n] \simeq & \textrm{Leray spectral sequence, $\lambda$ finite} \nonumber \\
R(q_{\widehat{Y}}\circ \lambda)_*(\mathcal{P})[n] \simeq & \nonumber \\
Rp_{\widehat{Y},*}(\mathcal{P})[n]. \phantom{\simeq}& \nonumber
\end{eqnarray}
This completes the proof of Corollary \ref{main cor} and the last part of Corollary \ref{corollary}. \hfill $\square$
\end{dimo1*}

\addcontentsline{toc}{chapter}{\bibname}
\printbibliography

\Addresses

\end{document}